\documentclass[11pt]{article}
\usepackage{amsmath,amssymb,amscd,amsthm}
\usepackage{epsfig}
\usepackage{graphicx}

\makeatletter\@addtoreset {equation}{section}\makeatother

\newtheorem{theorem}{Theorem}
\newtheorem{lemma}{Lemma}
\newtheorem{remark}{Remark}

\newtheorem{corollary}{Corollary}

\def\R{\mathbb{R}}

\def\N{\mathbb{N}}

\newenvironment{proof1}{
    \noindent {\it Proof }}{\hfill $\Box$}

\begin{document}

\title{\bf Long time dynamics and coherent states in nonlinear wave equations}

\author{E. Kirr\thanks{Department of Mathematics, University of Illinois,
Urbana--Champaign, Urbana, IL 61801, e-mail ekirr@illinois.edu}}

\date{April 29, 2016}

\maketitle

\begin{abstract}
We discuss recent progress in finding all coherent states supported by nonlinear wave equations, their stability and the long time behavior of nearby solutions.
\end{abstract}

\section{Introduction}\label{se:intro}
Since the discovery of wave equations two and a half centuries ago, many scientist and mathematicians have tried to understand their most striking feature, the coherent structures. An exact definition of coherent structures
will be given in Section \ref{se:gnw} but, formally, they are
solutions which propagate without changing shape (or with periodic change in shape). In the linear case, they are the eigenfunctions of the corresponding wave operator and, via spectral decomposition, the actual dynamics becomes a superposition of these coherent states (eigenfunctions) and a projection onto the remaining (continuous) spectrum. Once the latter had been analyzed via dispersive estimates or scattering wave operators, it became clear that the following asymptotic completeness conjecture is true for the linear case.
\begin{description}
\item[Asymptotic Completeness Conjecture:] Any initial data evolves towards a superposition of coherent structures plus a part that radiates (scatters) to infinity.
\end{description}
In the nonlinear case, some coherent structures are known either as minimizer or mountain pass type critical points of the energy subject to certain constraints, see Subsection \ref{sse:vm}. Other coherent states (sometimes the same ones) can be found via bifurcations from already known solutions such as the trivial one, see Subsection \ref{sse:bm}. In many situations the coherent states undergo symmetry breaking phenomena, see for example \cite{jw, kksw, AFGST:02, gus}, which are very important in practical applications. But none of these results, nor the mathematical methods they rely on can claim that they can actually identify all coherent states supported by a given nonlinear wave equation. Consequently, the asymptotic completeness conjecture is wide open with two notable exceptions: the case of completely integrable systems (such as the cubic Schr\"{o}dinger equation in one dimension) where a scattering transform renders the problem linear, or the case of weakly nonlinear regimes (i.e. small initial data) where at most two (small) coherent states are present and one of them is selected after a long transitional time, see \cite{ty1, ty2, ty3, swc}.

This paper aims to present recent results and new ideas for finding {\em all coherent states}
(solitary waves, breathers, kinks, vortices, etc) supported by a given nonlinear wave equations. As shown in Subsection \ref{sse:bm}, coherent states can be viewed as zeroes of a map between Banach spaces. Then, the global bifurcation theory, see \cite{bt, rab}, allows us to organize them in smooth
manifolds which either form loops or can be extended to the boundary of the domain inside which the linearization of the map is Fredholm. The new results and ideas concern finding all the limit points of such manifolds on the boundary of the Fredholm domain. Hence, from these limit points, the manifolds of coherent states {\em can be found and traced} both theoretically and numerically inside the Fredholm domain. Moreover, by comparing the spectrum
of the linearized operator at the two ``end points" of these manifolds we can deduce whether eigenvalues cross zero which is equivalent with the existence of bifurcations in this Fredholm region. Once discovered, the bifurcations can be studied using local bifurcation techniques to determine the new branches (manifolds) emerging from them and their dynamical stability. Global bifurcation theory now implies that the new branches have ``end points" on the boundary of the Fredholm domain. Consequently, we are able to find the bifurcations along the new branches and the process iterates until all these branches are discovered and matched with all possible limit points.

Results and open problems regarding the orbital stability of the manifolds of coherent states are discussed in Subsection \ref{sse:os} while Section \ref{se:as} discusses their asymptotic stability. The latter brings us closer to a resolution of the {\em Asymptotic Completeness Conjecture} but, unfortunately, it only describes the dynamics in a neighborhood of the coherent state manifolds. The last section is reserved for concluding remarks.

\section{General Hamiltonian Formulation}\label{se:gnw}

Most models related to wave propagation, in particular the
Schr\"{o}dinger, Hartree, Dirac, Klein-Gordon, Korteweg-de-Vries and
the classical wave equation, can be cast in the following general
framework, see \cite{gss}. The evolution of the quantity of interest $u$ is given
by:
 \begin{equation}\label{hpde}
 \frac{du}{dt}=JD_u{\cal E}(u),\qquad t\in\mathbb{R},
 \end{equation}
where $X$ is a real Hilbert space, the energy ${\cal E}:X\mapsto\mathbb{R}$ is a $C^2$
functional, $D_u$ denotes the Frechet derivative with respect to the
variable $u,$ and $J:D(J)\subseteq X^*\mapsto X$ is a skew-symmetric
operator $J^*=-J,$ defined on a dense subset of the dual of $X.$ Note that
even though $X$ is a Hilbert space, its dual is not necessarily
identified with $X.$ The reason is twofold: the applications have a
physically important larger Hilbert space $Y$ for which
$X\hookrightarrow Y=Y^*\hookrightarrow X^*,$ where all the
embeddings are dense, and the mathematical analysis of the operators
appearing in the applications rely on the larger Hilbert space $Y.$

Besides being time independent, the energy is in general invariant under additional groups of symmetries. These symmetries can be modeled by one or more (strongly) continuous groups of unitary operators. In what follows we will focus on one such group of symmetries $T(s):X\mapsto X,\ s\in\mathbb{R},$ which
leave the energy functional invariant and commutes with the $J$
operator:
$${\cal E}(T(s)u)={\cal E}(u),\qquad T(s)J=JT^*(-s),\quad\mbox{for all}\ s\in\mathbb{R},\ \mbox{and}\ u\in X.$$
By Noether's Theorem, see for example \cite{bf:bk},  the Hamiltonian dynamics has, besides energy, a second conserved quantity:
\begin{equation}\label{Qhpde}Q(u)=\frac{1}{2}\langle Bu,u\rangle,\qquad u\in X\end{equation}
provided that there exists a bounded, self-adjoint, linear operator
$B:X\mapsto X^*$ such that $JB$ extends the infinitesimal generator
of the continuous group: $T'(0).$

The coherent states are solutions of \eqref{hpde} of the type:
$$u(t)=T(\omega t)\phi_\omega$$
where $\phi_\omega\in D(T'(0))\subseteq X,\ \omega\in\mathbb{R}$ are
fixed. In applications $\phi_\omega$ usually gives the shape of $u$,
so these are solutions which do not change their shape as they
propagate. By plugging in \eqref{hpde} one finds that:
$$JD_\phi {\cal E}(\phi_\omega)=J\omega D_\phi Q(\phi_\omega).$$
Consequently the solutions in $D(T'(0))$ of the stationary equation:
\begin{equation}\label{cseq}
D_\phi {\cal E}(\phi)=\omega D_\phi Q(\phi)
\end{equation}
generate coherent states and they are the only possible coherent structures if $J$ is one-to-one.

Coherent states are orbitally stable if any solution of
\eqref{hpde} starting close to the orbit $T(s)\phi_\omega,\
s\in\mathbb{R},$ of the coherent state, remains close to it at all
times. More precisely for any $\varepsilon$ there exists a $\delta$
such that
$$\inf_{s\in\mathbb{R}}\|u(0)-T(s)\phi_\omega\|<\delta\quad\mbox{implies}\quad\sup_{t\in\R}\inf_{s\in\R}\|u(t)-T(s)\phi_\omega\|<\varepsilon.$$

Asymptotic stability means certain convergence of the solutions to
the orbit of a coherent structure and usually takes the form: there
exists a Banach space $Z,$ $X\hookrightarrow Z$ densely, and
$\delta>0$ such that if
$\inf_{s\in\mathbb{R}}\|u(0)-T(s)\phi_\omega\|_{X\bigcap
Z^*}<\delta$ then there exists a coherent structure $T(\omega_+
t)\phi_{\omega_+}$ (close to $T(\omega t)\phi_\omega$) with the
property:
$$\lim_{t\rightarrow\infty}\inf_{s\in\R}\|u(t)-T(s)\phi_{\omega_+}\|_{Z}=0.$$

For example, in the case of the nonlinear Schr\" odinger equation
(NLS) we have $X=H^1(\mathbb{R}^n),$ the Sobolev space of complex
valued functions but with the real Hilbert space structure,
$X^*=H^{-1},\ Jv=-iv,\ T(s)u=e^{-is}u,$
\begin{eqnarray}
{\cal E}(u)&=&\frac{1}{2}\int_{\mathbb{R}^n}|\nabla u(x)|^2dx+\frac{1}{2}\int_{\mathbb{R}^n}V(x)|u(x)|^2dx+\frac{\gamma}{p+2}\int_{\mathbb{R}^n}|u(x)|^{p+2}dx,\label{Enls}\\
Q(u)&=&\frac{1}{2}\int_{\mathbb{R}^n}|u(x)|^2dx,\label{Qnls}
\end{eqnarray}
hence the evolution equation \eqref{hpde} becomes:
\begin{equation}\label{eq:nse}
i\frac{\partial u}{\partial t}=(-\Delta + V(x))u(t,x)+\gamma
|u|^{p}u(t,x),\quad t\in\R,\ x\in\R^n,\end{equation}
while the coherent structures are solutions of the form $u(t,x)=e^{iEt}\phi_E(x),\ E=-\omega\in\R,\ \phi_E\in H^1(\R^n),$ and satisfy the equation:
\begin{equation}\label{eq:ev}
F(\phi_E,E)=(-\Delta + V+E)\phi_E+\gamma
|\phi_E|^{p}\phi_E=0.
\end{equation}
Here $V:\R^n\mapsto\R$ is called the potential, $\gamma\in\R$ measures the strength of the nonlinear interaction while its sign classifies it into attractive for $\gamma <0,$ and repelling for $\gamma>0.$ 
In this context the coherent states are usually called bound states or, in the translation invariant case $V\equiv 0,$ solitons. The Hartree Equation has exactly the same $X,\ J,\ T$ and $Q$ but the superquadratic term in the energy becomes nonlocal:
\begin{equation}\label{Ehar}
{\cal E}(u)=\frac{1}{2}\int_{\mathbb{R}^n}|\nabla u(x)|^2dx+\frac{1}{2}\int_{\mathbb{R}^n}V(x)|u(x)|^2dx+\frac{\gamma}{4}\int_{\mathbb{R}^n}K(x,y)|u(x)|^2|u(y)|^2dxdy,
\end{equation}
where the kernel $K\geq 0.$

\section{Coherent States}
This section illustrates how one can find all coherent states of equations of type \eqref{hpde} i.e., all solutions of \eqref{cseq}, and how one can determine their orbital stability. Traditionally, large coherent structures are found via variational methods, for example as minimizers of the energy subject to a fixed value of the second conserved quantity. However, as we shall see in the next subsection, the variational techniques are not capable of finding all coherent states. Instead we will show in subsection \ref{sse:bm} how bifurcation methods, in particular the analytical global bifurcation theory \cite{bt}, can be enhanced to determine all coherent states, and their orbital stability, see subsection \ref{sse:os}.

 \subsection{Variational Methods. Existence and Stability of Ground States.}\label{sse:vm}
The coherent states equation \eqref{cseq} coincides with the equation for the critical points of the energy ${\cal E}$ restricted to the level sets of the second conserved quantity $Q$ i.e., it is the Euler-Lagrange equation for the energy subject to the constrain $Q=constant.$ An important subset of the coherent states are the ground states which are minimizers of the energy under the constrain:
\begin{equation}\label{eq:var}{\cal E}(\phi_\omega)=\min_{\phi\in X,Q(\phi)=\mu}{\cal E}(\phi),\quad\mbox{and}\quad \mu\in\R \end{equation}
If the energy (subject to the constrain) is bounded from below:
$$\exists m\in\R\mbox{ such that }{\cal E}(\phi)\geq m\ \forall\phi\in X\mbox{ with }Q(\phi)=\mu,$$
and coercive: $$ \lim_{\|\phi\|_X\rightarrow\infty,Q(\phi)=\mu}{\cal E}(\phi)=\infty,$$ then minimizing sequences $\{\phi_n\}_{n\in\N}\subset X$ are bounded and, due to the reflexivity of the Hilbert space $X,$ have at least one weak limit, say $\phi_{n_k}\rightharpoonup\phi_0.$ For $\phi_0$ to be a ground state it must satisfy
\begin{equation}\label{gscond}Q(\phi_0)=\mu=\lim_{k\rightarrow\infty}Q(\phi_{n_k}),\quad {\cal E}(\phi_0)\leq\lim_{k\rightarrow\infty}{\cal E}(\phi_{n_k})=\inf_{\phi\in X,Q(\phi)=\mu}{\cal E}(\phi)\end{equation}
Note that these conditions are not trivially satisfied as $Q$ or ${\cal E}$ might not be weakly sequentially continuous even though they are both continuous with respect to the norm on $X.$ These two issues are resolved by compactness arguments which show that weakly convergent minimizing sequences are actually strongly convergent i.e.,
\begin{equation}\label{compcd}\phi_{n_k}\stackrel{X}{\rightharpoonup}\phi_0\quad\mbox{ implies }\quad\lim_{k\rightarrow\infty}\|\phi_{n_k}-\phi_0\|_X=0,\end{equation}
see for example \cite{bf:bk,caz:bk}.

In NLS with confining potentials $V\geq 0,\ \lim_{|x|\rightarrow\infty}V(x)=\infty,$ the potential restricts the domain of finite energy to a compact subspace. More precisely we have:
\begin{equation}\label{Xnlsconf}
X=\left\{\phi\in H^1(\R^n)\ |\ \int_{\mathbb{R}^n}V(x)|\phi(x)|^2dx\right\}\mbox{ with } \langle\phi,\psi\rangle_X=\langle\phi,\psi\rangle_{H^1}+\int_{\mathbb{R}^n}V(x)\overline\phi(x)\psi(x)dx,\end{equation}
is a Hilbert space which satisfies:
\begin{equation}\label{Xnlscomp} X\stackrel{compact}{\hookrightarrow}L^p(\R^n),\ 2\leq p<\frac{2n}{n-2}\mbox{ if }n\geq 3,\ 2\leq p<\infty\mbox{ if } n=1,2.\end{equation}
In particular, the above weakly convergent, minimizing subsequence:
$$\phi_{n_k}\stackrel{X}{\rightharpoonup}\phi_0\mbox{ is strongly convergent in }L^p(\R^n)\mbox{ i.e., } \lim_{k\rightarrow\infty}\|\phi_{n_k}-\phi_0\|_{L^p}=0,\mbox{ for }2\leq p<\frac{2n}{n-2}.$$
Therefore $$\mu=Q(\phi_{n_k})=\frac{1}{2}\int_{\Omega}|\phi_{n_k}(x)|^2dx\stackrel{k\rightarrow\infty}{\longrightarrow}\frac{1}{2}\int_{\Omega}|\phi_0(x)|^2dx=Q(\phi_0)$$
and $$\int_{\Omega}|\phi_{n_k}(x)|^{p+2}dx\stackrel{k\rightarrow\infty}{\longrightarrow}\int_{\Omega}|\phi_0(x)|^{p+2}dx.$$ Combining the above with the weak lower semicontinuity of the first two (kinetic and potential) terms in the energy which are convex, we deduce that \eqref{gscond} holds and $\phi_0$ is a ground state. Moreover, the following inequality holds
$${\cal E}(\phi_0)\geq \inf_{\phi\in X,Q(\phi)=\mu}{\cal E}(\phi)=\lim_{k\rightarrow\infty}{\cal E}(\phi_{n_k}),$$
which is opposite to \eqref{gscond}. Therefore, on this minimizing subsequence, the kinetic and potential terms must be convergent which combined with the convergence of $Q$ gives $\|\phi_{n_k}\|_X\rightarrow\|\phi_0\|_X$ in addition to $\phi_{n_k}\rightharpoonup\phi_0.$ The strong convergence \eqref{compcd} now follows from the uniform convexity of the Hilbert space $X.$

Essential in the above argument is the compactness of the embeddings \eqref{Xnlscomp}. Heuristically, one might think that $\int_{\mathbb{R}^n}V(x)|\phi(x)|^2dx < \infty$ and $\lim_{|x|\rightarrow\infty}V(x)=\infty$ forces a ``uniform decay at infinity" on $\phi\in X$ which does imply compactness, see \cite[Section 1.7]{caz:bk}. But this is not quite correct since if $\phi$ is a countable sum of of smoothed characteristic functions of disjoint annuli with the same center and exterior and interior radius growing to infinity we have $\lim_{|x|\rightarrow\infty}\phi(x)\not=0,$ but $\phi\in X$ if the Lebesgue measure (volume) of the annuli converges to zero sufficiently fast. However, the argument in \cite[Section 1.7]{caz:bk} can be adapted to prove \eqref{Xnlscomp} as follows. Consider an arbitrary bounded sequence $$\{\phi_n\}_{n\in\N}\subset X\mbox{ with }\|\phi_n\|_X\leq M,\mbox{ for all }n\in\N$$ Since $X$ is Hilbert the sequence has a weakly convergent subsequence $$\phi_{n_k}\stackrel{X}{\rightharpoonup}\phi_0\in X.$$ Then, for each $\varepsilon>0$ we can choose  $R>0$ such that $$V(x)> 16\varepsilon^{-2}\max\left\{M^2,\int_{\mathbb{R}^n}V(x)|\phi_0(x)|^2dx\right\},\mbox{ for }|x|> R.$$ Consequently, we have
$$M^2\geq\int_{\mathbb{R}^n}V(x)|\phi_{n_k}(x)|^2dx\geq\int_{|x|> R}V(x)|\phi_{n_k}(x)|^2dx\geq \frac{16M^2}{\varepsilon^2}\int_{|x|> R}|\phi_{n_k}(x)|^2dx,$$ and
$$\int_{\mathbb{R}^n}V(x)\|\phi_{0}(x)\|^2dx\geq\int_{|x|> R}V(x)|\phi_{0}(x)|^2dx\geq \frac{16\int_{\mathbb{R}^n}V(x)|\phi_{0}(x)|^2dx}{\varepsilon^2}\int_{|x|> R}|\phi_{0}(x)|^2dx,$$
which imply
 $$\left(\int_{|x|> R}|\phi_{n_k}(x)-\phi_0(x)|^2dx\right)^{1/2}\leq\left(\int_{|x|> R}|\phi_{n_k}(x)|^2dx\right)^{1/2}+\left(\int_{|x|> R}|\phi_0(x)|^2dx\right)^{1/2}\leq \varepsilon/2.$$ Now:
\begin{eqnarray}
\|\phi_{n_k}-\phi_0\|_{L^2(\mathbb{R}^n)}&=&\|\phi_{n_k}-\phi_0\|_{L^2(|x|<R)}+\|\phi_{n_k}-\phi_0\|_{L^2(|x|>R)}<\|\phi_{n_k}-\phi_0\|_{L^2(|x|<R)}+\varepsilon/2.\nonumber
\end{eqnarray}
But, by Rellich-Kondrachov Theorem, $H^1(|x|<R)$ is compactly embedded in $L^2(|x|<R)$ which means that the $X$ hence $H^1$ weakly convergent sequence $\{\phi_{n_k}\}$ is strongly convergent in $L^2(|x|<R)$ and we can choose $k(\varepsilon)\in\N$ such that $\|\phi_{n_k}-\phi_0\|_{L^2(|x|<R}<\varepsilon/2$ for all $k>k(\varepsilon).$ All in all, for each $\varepsilon>0$ we can find $k(\varepsilon)\in\N$ such that $\|\phi_{n_k}-\phi_0\|_{L^2(\R^n)}<\varepsilon$ for all $k>k(\varepsilon)$ i.e., $\phi_{n_k}$ converges strongly (in norm) in $L^2(\R^n).$ Moreover, $\phi_{n_k}$ bounded in $X$ hence in $H^1(\R^n)$ also implies, via Sobolev embedding, that it is bounded in $L^{2n/(n-2)}(\R^n)$ and, by interpolation, convergent to $\phi_0$  in $L^p,\ 2\leq p<2n/(n-2).$ So, any bounded sequence in $X$ has a convergent subsequence in $L^p(\R^n),\ 2\leq p<2n/(n-2)$ if $n\geq 3,\ 2\leq p <\infty$ if $n=1,2,$ which implies \eqref{Xnlscomp}.

However, in general, the verification of \eqref{gscond} requires {\em concentration compactness}, see \cite{PLio1,PLio2} or \cite[Section 1.7]{caz:bk}. This theory will be discussed in a different context in the next subsection. Suffices to say that in the NLS example it covers the case of non-confining potentials: $\lim_{|x|\rightarrow\infty}V(x)=0,$ (when the energy is bounded from below.)

The main difference between the ground states given by \eqref{eq:var} and other solutions of \eqref{cseq} (called excited states) is that the former are in general stable under the dynamics:

\begin{theorem}\label{th:stabgs} Fix $\mu\in\R$ and assume the set of ground states, $G=\{\phi_\omega\in X\ |\ \phi_\omega\mbox{ solves \eqref{eq:var}}\},$ is non-empty. Fix $\phi_0\in G$ and further assume that any minimizing sequence of \eqref{eq:var} has a strongly convergent subsequence in $X.$ Then for any $\varepsilon>0$ there exists $\delta>0$ such that for all $u_0\in X$ with $\|u_0-\phi_0\|_X<\delta$ we have that the solution $u(t)$ of the wave equation \eqref{hpde} with initial condition $u(0)=u_0$ remains within $\varepsilon$ distance from $G$ for all times.
\end{theorem}

\begin{proof1} Suppose contrary, there is an $\varepsilon>0,$ a sequence $\{u_n\}_{n\in\N}\subset X$ with $\|u_n-\phi_0\|_X\rightarrow 0$ and a sequence of times $\{t_n\}_{n\in\N}\subset\R$ such that the solutions $u_n(t)$ of the wave equation \eqref{hpde} with initial condition $u(0)=u_n$ satisfy
$${\rm dist}(u_n(t_n),G)=\inf\{\|u(t_n)-\psi\|_X\ |\ \psi\in G\} \geq\varepsilon.$$
By continuity of $Q,$ we have $Q(u_n)\rightarrow Q(\phi_0)=\mu,$ and, by using its bilinear form \eqref{Qhpde}, we can find $\{\lambda_n\}_{n\in\N}\subset\R$ such that
$$Q(\lambda_nu_n)=\lambda^2_nQ(u_n)=\mu\quad\mbox{and}\quad\lambda_n\rightarrow 1.$$
We now claim that $\{\lambda_nu_n(t_n)\}_{n\in\N}\subset X$ is a minimizing sequence for \eqref{eq:var}. Indeed, by conservation of $Q$ along solutions of \eqref{hpde} we have $Q(u_n(t))=Q(u_n)$ for all $t\in\R$ in particular:
$$Q(\lambda_nu_n(t_n))=\lambda^2_nQ(u_n(t_n))=\lambda^2_nQ(u_n)=\mu,$$
while by conservation of the energy we have:
$${\cal E}(\lambda_nu_n(t_n))={\cal E}(\lambda_nu_n)\rightarrow {\cal E}(\phi_0)=\min_{\phi\in X,Q(\phi)=\mu}{\cal E}(\phi),$$
where the convergence follows from the continuity of ${\cal E}$ w.r.t. the norm in $X.$

Now, since $\{\lambda_nu_n(t_n)\}_{n\in\N}\subset X$ is a minimizing sequence for \eqref{eq:var}, it has, by hypothesis, a convergent subsequence to some $\phi_1\in X$ i.e., $\|\lambda_nu_n(t_n)-\phi_1\|_X\rightarrow 0.$ But, by continuity of $Q$ and ${\cal E}$ we have:
$$Q(\phi_1)=\lim_{k\rightarrow\infty}Q(\lambda_{n_k}u_{n_k}(t_{n_k}))=\mu,\ {\cal E}(\phi_1)=\lim_{k\rightarrow\infty}{\cal E}(\lambda_{n_k}u_{n_k}(t_{n_k}))={\cal E}(\phi_0)=\min_{\phi\in X,Q(\phi)=\mu}{\cal E}(\phi)$$
i.e. $\phi_1\in G.$ Therefore we have:
$${\rm dist}(u_n(t_n),G)\leq \|u_n(t_n)-\phi_1\|_X\leq |1-\lambda_n|\|u_n\|_X+\|\lambda_nu_n(t_n)-\phi_1\|_X\rightarrow 0,$$
since $\lambda_n\rightarrow 1$ and $\|\lambda_nu_n(t_n)-\phi_1\|_X\rightarrow 0.$ This contradicts our assumption that ${\rm dist}(u_n(t_n),G)\geq\varepsilon$ and finishes the proof of the theorem.
\end{proof1}

\begin{remark} Note that all examples discussed above satisfy the hypotheses of the Theorem. Moreover, with the exception of the case $V\equiv 0$ the set of ground states is unique up to the symmetries induced by the semigroup $T:$
\begin{equation}\label{gsuniq} G=\{T(s)\phi_0\ |\ \mbox{for some }\phi_0\mbox{ which solves \eqref{eq:var} and all }s\in\R\},\end{equation}
provided $\mu$ is small, see for example \cite{RW:88,AFGST:02}.
\end{remark}
Note that the invariance of both $Q$ and ${\cal E}$ w.r.t $T$ automatically implies $G\supseteq\{T(s)\phi_0\ :\ s\in\R\}$ if $\phi_0$ solves \eqref{eq:var}. However the equality between the two sets implies orbital stability:

\begin{corollary}\label{cor:stabgs}
Under the assumptions in Theorem \ref{th:stabgs}, if, in addition, \eqref{gsuniq} holds, then the ground state $\phi_0$ is orbitally stable.
\end{corollary}
\begin{proof1} By the theorem, for each $\varepsilon>0$ there exists $\delta>0$ such that for all $u_0\in X$ with $\|u_0-\phi_0\|_X<\delta$ we have:
$$\sup_{t\in\R}{\rm dist}(u(t),G)<\varepsilon.$$
But, in this case
$${\rm dist}(u(t),G)=\inf_{s\in\R}\|u(t)-T(s)\phi_0\|_X$$ which implies orbital stability, see the definition below \eqref{cseq}.
\end{proof1}

\begin{remark}More generally, the orbital stability result in the previous corollary holds even if the ground states are not unique (up to the action of $T$) provided that the orbit $T(s)\phi_0$ is separated from the orbits of the other ground states by a fixed distance $d>0.$ Just use $\varepsilon <d/2$ in the above proof and note that the set of points at distance less the $d/2$ from the orbit of $\phi_0$ is disjoint from the set of all points at distance less the $d/2$ from the orbits of all other ground states while the map ${\rm dist}(u(t),G)$ is continuous in time.
\end{remark}

In the case $V\equiv 0,$ it is necessary to mode out the second (hidden) symmetry, namely the invariance of both energy and $Q$ with respect to translations, see for example \cite[Section 8.3]{caz:bk}. One obtains:
$$\sup_{t\in\R}\inf_{s\in\R,y\in\R^n}\|u(t)-T(s)\phi_0(\cdot+y)\|<\varepsilon.$$

Another advantage of the global minimization problem \eqref{eq:var} over other solutions of \eqref{cseq} is that certain manipulation of the functions $\phi\in X,$ such as taking its absolute value or symmetrizing it, can lower the energy and provide information on the shape of the ground states. For example, in NLS the ground states are always positive $\phi_\omega>0$ up to the action of the group $T$ i.e., up to rotations in the complex plane, and, in the absence of potential $V\equiv 0,$ they are spherically symmetric,
see \cite[Chapter 8]{caz:bk}. Moreover, in \cite{AFGST:02}, see also \cite{gus} for a related result, the authors show that the ground states for the Hartree equation must localize at global minima of the potential in the limit $\mu\rightarrow\infty,$ and, consequently, the following symmetry breaking phenomena occurs:

\begin{theorem}\label{th:sb1} Consider the Hartree example \eqref{Ehar} with an attractive nonlinearity $\gamma <0,$ and a continuous, bounded potential $V$ which is invariant under a finite group of Euclidian symmetries on $\R^n.$ If the action of the group is nontrivial on any global minima of $V$ then the are $\mu_0\leq \mu_1$ such that the ground states with $\mu<\mu_0$ are invariant under the group of Euclidian symmetries but the ground states with $\mu>\mu_1$ are not invariant.
\end{theorem}

The disadvantages of the minimization problem \eqref{eq:var} are that it cannot give all coherent states i.e., all solutions of \eqref{cseq}, and it may have no solutions. This is the case when the energy is not be bounded from below (even when constrained to $Q=const$) which occurs in the NLS example with critical and supercritical nonlinearities $p\geq 4/n.$ The issue is sometimes resolved by reformulating the problem i.e., by finding the global minimum of a functional different from the energy. In \cite{RW:88} the authors use the following reformulation:
$$\min_{\phi\in X,\phi\not=0}J^E(\phi)=\frac{\int_{\R^n}|\nabla\phi(x)|^2+V(x)|\phi(x)|^2+E|\phi(x)|^2}{\left(\int_{\R^n}|\phi(x)|^{p+2}dx\right)^{\frac{2}{p+2}}}$$
which is equivalent to:
$$\min_{\phi\in X,\int_{\R^n}|\phi(x)|^{p+2}dx=\mu}\int_{\R^n}|\nabla\phi(x)|^2+V(x)|\phi(x)|^2+E|\phi(x)|^2,\qquad \mu\in\R ,$$
see \cite[Chapter 8]{caz:bk} for other possible reformulations. An important result in \cite{RW:88} is the existence of a solution to \eqref{eq:ev} for each $E>E_0$ where $-E_0$ is the lowest eigenvalue of $-\Delta+V.$ However, the method is still limited to finding a subset of solutions of \eqref{eq:ev} and does not provide any direct information regarding their dynamical stability (the latter may be fixed by combining information on the Hessian of $J^E$ at a minima with the techniques described in subsection \ref{sse:os}). One expects that these variational reformulations will also lead to information on the shape and localization of the ground states, consequently symmetry breaking results may be proven for critical and supercritical nonlinearities, see \cite{gus}. While other critical points of the energy or other associated functionals may be found via mountain pass techniques, see for example \cite{ny}, the variational methods do not provide a systematic method to identify all solutions of \eqref{cseq}.

 \subsection{Bifurcation Methods}\label{sse:bm}

This section discusses recent progress towards finding all coherent states of a nonlinear wave equation i.e., all solutions of equation \eqref{cseq}. We will focus on the NLS example for which \eqref{cseq} becomes \eqref{eq:ev}:
$$
F(\phi_E,E)=(-\Delta + V+E)\phi_E+\gamma
|\phi_E|^{p}\phi_E=0,
$$
where the potential $V:\R^n\mapsto\R,$ is first assumed to be non-confining, $\lim_{|x|\rightarrow\infty}V(x)=0,\ V\in L^q(\R^n)+L^\infty(\R^n)$ for some $q\geq 1,\ q>n/2.$ The case of confining potentials is discussed in Remark \ref{rmk:pot}. We will assume that the power of the nonlinearity satisfies $0<p<\infty,$ if $n=1$ or 2 and $0<p<4/(n-2)$ if $n\geq 3,$ which insures local well posedness of the time dependent equation \eqref{hpde} with initial data in the Sobolev space $H^1(\R^n).$ Of special interest are the ground states which, for the purpose of this presentation, will be defined as coherent states i.e., solutions of \eqref{eq:ev} that satisfy $\phi_E(x)>0,\ \forall x\in\R^n,$ modulo multiplication by a complex number of modulus one (modulo rotations). Note that, in NLS this definition encompasses the ground states discussed in the previous subsection and characterized by \eqref{eq:var} because the latter satisfy the positivity condition. Now all coherent states can be viewed as zeroes of the map
$F:H^1(\R^n)\times\R\mapsto H^{-1}(\R^n)$ where the Sobolev spaces
are endowed with their {\em real} Hilbert space structure in order
for $F$ to be differentiable. Note that $F$ is equivariant under rotations: $$F(e^{i\theta}\phi,E)=e^{i\theta}F(\phi,E),\ \theta\in\R,$$ hence the solution set for \eqref{eq:ev} is invariant under rotations.

We will study the solutions of \eqref{eq:ev} in the subdomain of $H^1(\R^n)\times\R$ where its linearization is Fredholm. We have
\begin{eqnarray}\lefteqn{D_\phi F(\phi,E)[u+iv]=}\nonumber\\
&&\left[\begin{array}{lr} -\Delta+V+E+\gamma(p+1)|\phi|^p-2\gamma(\Im\phi)^2|\phi|^{p-2} & 2\gamma(\Re\phi)(\Im\phi)|\phi|^{p-2}\\
2\gamma(\Re\phi)(\Im\phi)|\phi|^{p-2} & -\Delta+V+E+\gamma|\phi|^p+2\gamma(\Im\phi)^2|\phi|^{p-2}\end{array}\right]\left[\begin{array}{c} u\\ v\end{array}\right]\nonumber\end{eqnarray}
where we separated the real and imaginary parts of the complex valued functions involved into the first and second component. So,
 $$D_\phi F(\phi,E)=\left[\begin{array}{lr} -\Delta+E & 0\\
0 & -\Delta+E \end{array}\right]+{\cal V}(\phi)$$
where, for any $\phi\in H^1(\R^n),$
$${\cal V}=\left[\begin{array}{lr} V+\gamma(p+1)|\phi|^p-2\gamma(\Im\phi)^2|\phi|^{p-2} & 2\gamma(\Re\phi)(\Im\phi)|\phi|^{p-2}\\
2\gamma(\Re\phi)(\Im\phi)|\phi|^{p-2} & V+\gamma|\phi|^p+2\gamma(\Im\phi)^2|\phi|^{p-2}\end{array}\right].$$
is a relatively compact perturbation of the diagonal operator $-\Delta +E$ on $H^{-1}\times H^{-1}$ with domain $H^1\times H^1$ (or on $L^2\times L^2$ with domain $H^2\times H^2$), see for example \cite{RS:bk4}. Since the latter has essential (continuous) spectrum the interval on the real line  $[E,\infty)$ we get via Weyl's Theorem:
\begin{lemma}\label{th:R1}
$D_\phi F(\phi, E)$ is Fredholm (of index zero) iff $E>0.$ At the left boundary, $E=0,$ zero is at the edge of its essential (continuous) spectrum while for $E<0$ zero is inside the essential (continuous) spectrum.
\end{lemma}
We will restrict ourselves to the domain $H^1\times (0,\infty),$ i.e. $\{E>0\},$ which will now be called the bifurcation diagram. Note that, for $E<0,$ by the limiting absorbtion principle, there are no nontrivial solutions of \eqref{eq:ev} under mild assumptions on their decay rates, see \cite{bl}. Obviously, $(\phi_E\equiv 0,E),\ E\in\R$ solves \eqref{eq:ev}.

\begin{figure}
\begin{center}
\includegraphics[width=3in,height=6cm]{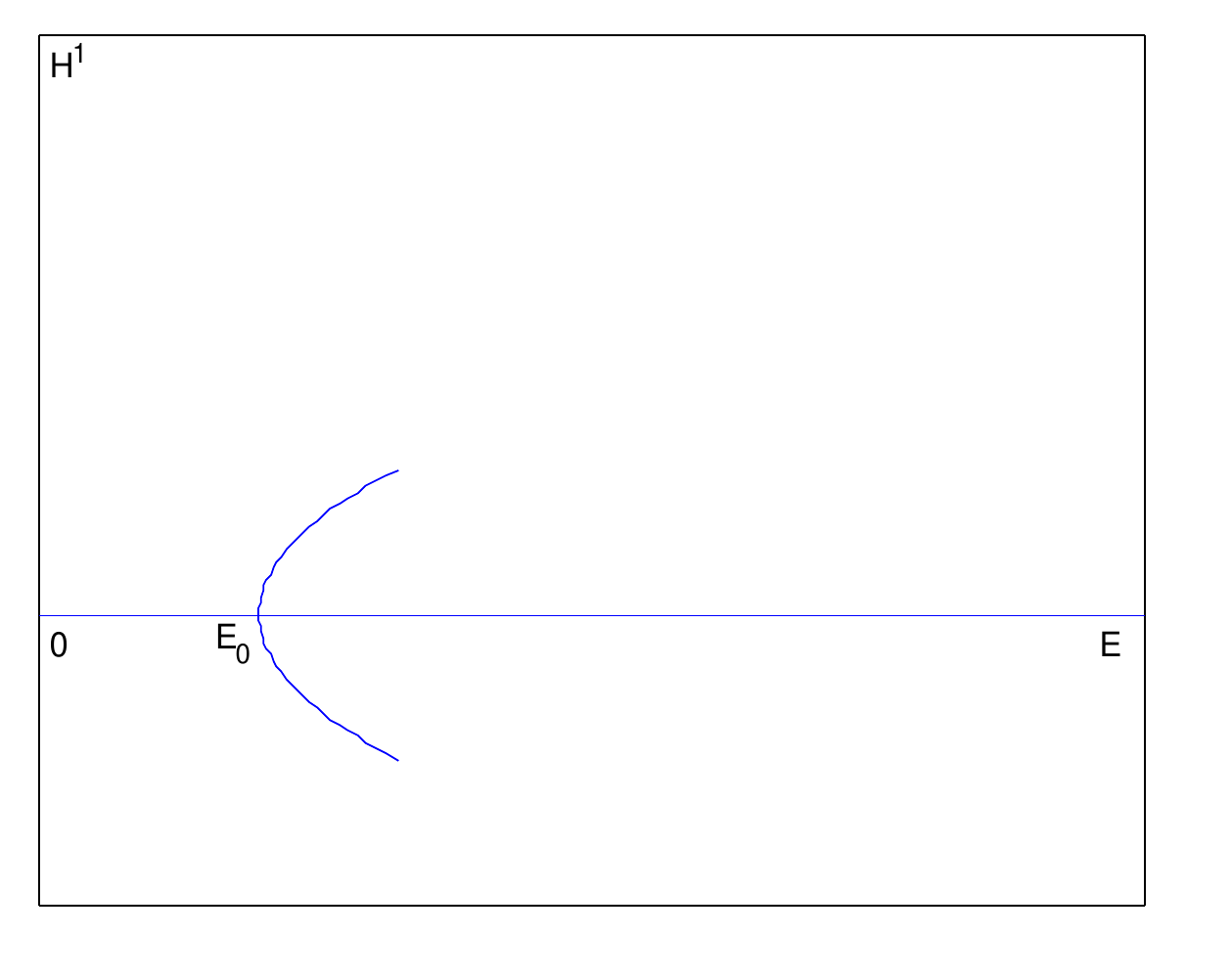}
\end{center}
\caption{\small Bifurcation diagram for bound states. Only the trivial zero
coherent state and the small ground states are represented for the attractive case $\gamma<0.$ For the repelling case, $\gamma>0,$ the branch points the other way, i.e. $E<E_0$ along the branch.}
\label{figbd}
\end{figure}

\begin{figure}
\begin{center}
\includegraphics[width=4in,height=10cm]{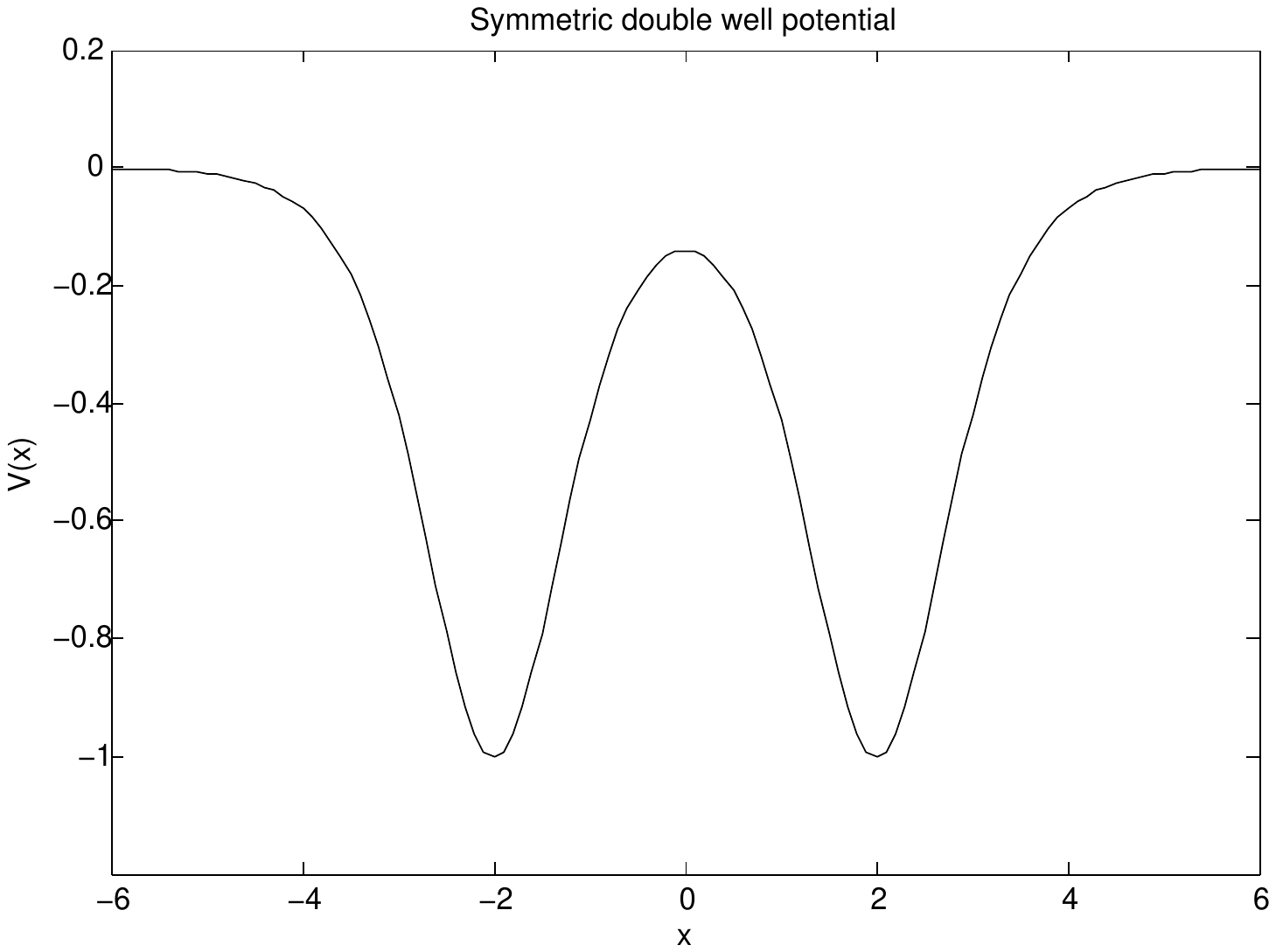}
\includegraphics[width=6in,height=6cm]{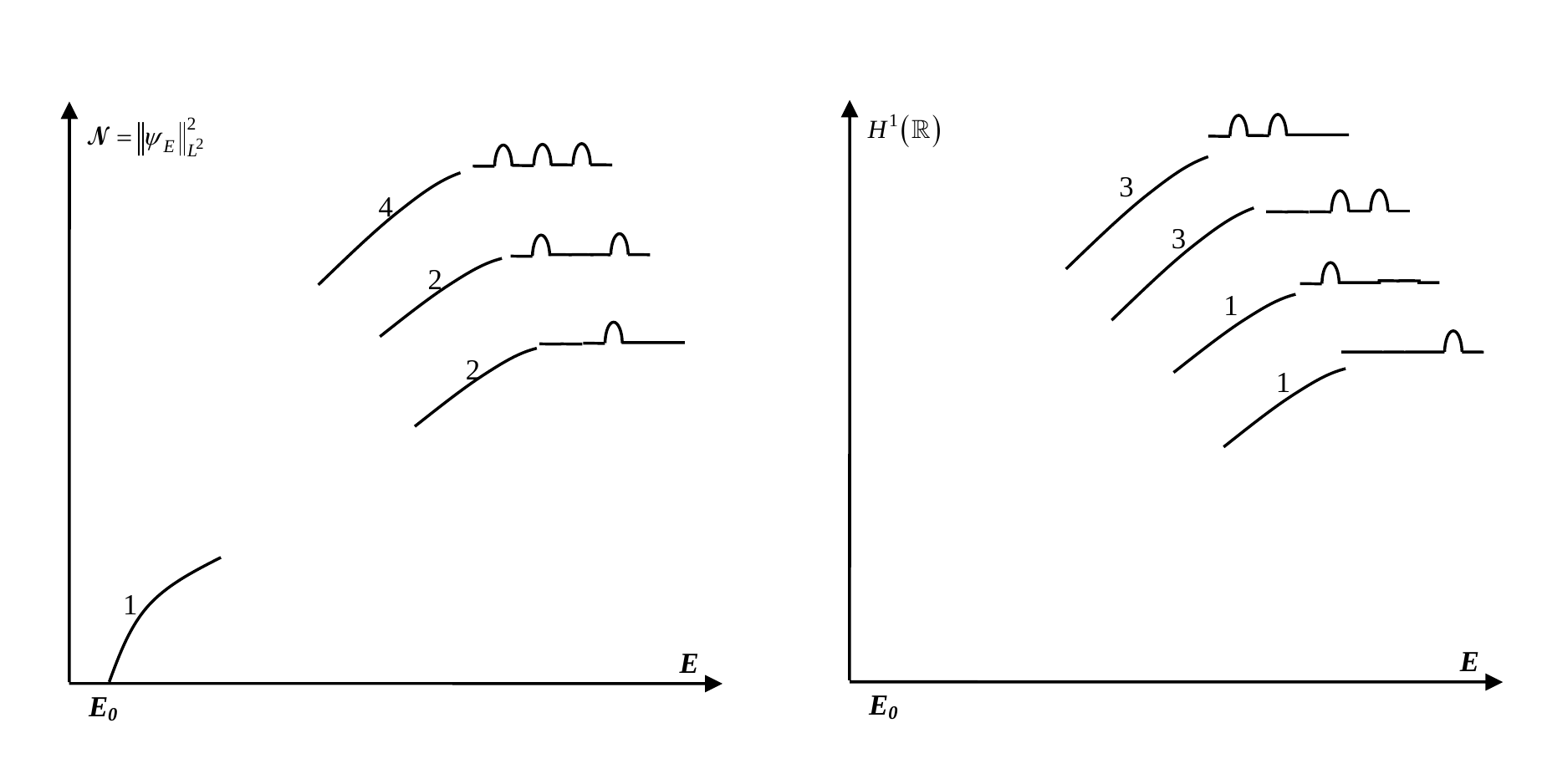}
\end{center}
\caption{An one dimensional even potential (top panel) and a sketch
of the corresponding symmetric (bottom left panel) and asymmetric
(bottom right panel) ground state branches for subcritical
nonlinearity, $p<2.$ The number on top of each branch gives the
number of negative eigenvalues for the linearized operator while the shape on the
right shows the actual shape of the solution on the branch.}
\label{figbranches}
\end{figure}

For a while we will assume:
\begin{description}
\item[(SA)] $-\Delta+V$ has at least one negative eigenvalue. The lowest will be denoted by $-E_0.$
\end{description}
Note that the assumption holds in space dimensions $n=1,2$ for non-trivial, negative potentials and requires potentials with sufficiently large negative parts in dimensions $n\geq 3.$ As shown for example in \cite{PW,kz1,kksw}, hypothesis (SA) leads to a pitchfork bifurcation at $(\phi_{E_0}\equiv 0,E_0),$ which creates exactly one curve (modulo rotations) of non-trivial ground states, see Figure \ref{figbd}. 
Moreover, if $V$ is invariant under a group of symmetries then so are the ground states on this branch. In particular, if $V$ is a symmetric, double well potential, see Figure \ref{figbranches} top panel, then the profile of the ground states is equally distributed between the two wells.

We are going to rely on global bifurcation theory which, besides a Fredholm linearization, also requires compactness either of the solution set of \eqref{eq:ev} (in the analytic case, see \cite{bt}) or of the map $\tilde F$ (in the continuous case where degree theory is used, see \cite{rab}) where $\tilde F$ is obtained by transforming \eqref{eq:ev} into a fixed point problem, for example:
\begin{equation}\label{eq:tF}\phi=(-\Delta + 1)^{-1}\psi,\ \psi=(1-E)\phi-V(x)\phi-\gamma
|\phi|^{p}\phi\stackrel{def}{=}\tilde F(\psi, E),\ \tilde F:H^{-1}\times\R\mapsto H^{-1}\end{equation}
Note that for $\phi$ defined on a bounded domain, $\Omega\subset\R^n,$ compactness of $\tilde F$ follows from compactness of Sobolev embeddings $H^1_0(\Omega)\hookrightarrow H^{-1}(\Omega),\ L^q(\Omega)\hookrightarrow H^{-1}(\Omega),\ 2\leq q<2n/(n-2),$ however on $\R^n$ the situation is much more delicate. For repelling nonlinearities, $\gamma >0,$ the problem is analyzed in \cite{jls:bsd} where the authors prove uniform bounds on the solutions of the inequality $|F(\phi(x),E)|\leq \Psi(x)$ to obtain the compactness necessary for defining a degree for $\tilde F$ given in \eqref{eq:tF}. Then global bifurcation theory implies that  from $(0,E),$ where $-E$ is any negative eigenvalue of $-\Delta +V$ with odd multiplicity, bifurcate branches of non-trivial solutions of \eqref{eq:ev} which``end up" either at the boundary of the bifurcation diagram, or at $(0,E_1)$ where $-E_1$ is a different eigenvalue of $-\Delta+V.$ Note that the results cannot give any information on existence of other bifurcation points along this branches, or on existence of other branches that may not connect to the trivial solution, or on the exact region or point where each branch ends up. In the particular case of the ground states bifurcating from $(0,E_0),$ the authors of \cite{jls:gsd} use further comparison theorems to infer that the branch approaches the left boundary $\{E=0\}$ of the bifurcation diagram, i.e. the part of the boundary where zero is at the edge of the continuous spectrum of the linearized operator, see Lemma \ref{th:R1} and Figure \ref{figbd} but note that for $\gamma>0$ the branch points the other way.

\begin{figure}
\begin{center}
\includegraphics[width=6in,height=6cm]{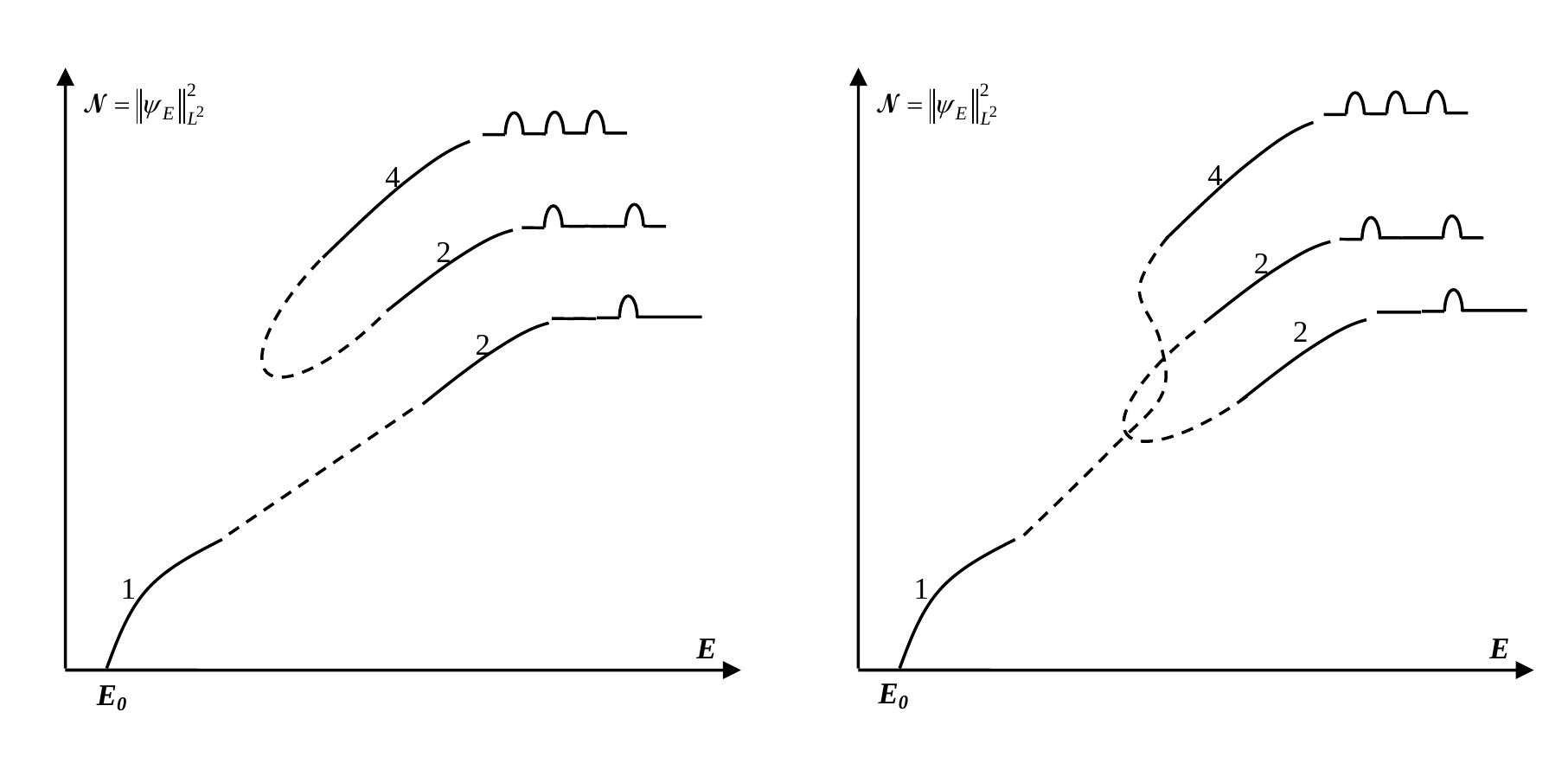}
\includegraphics[width=6in,height=6cm]{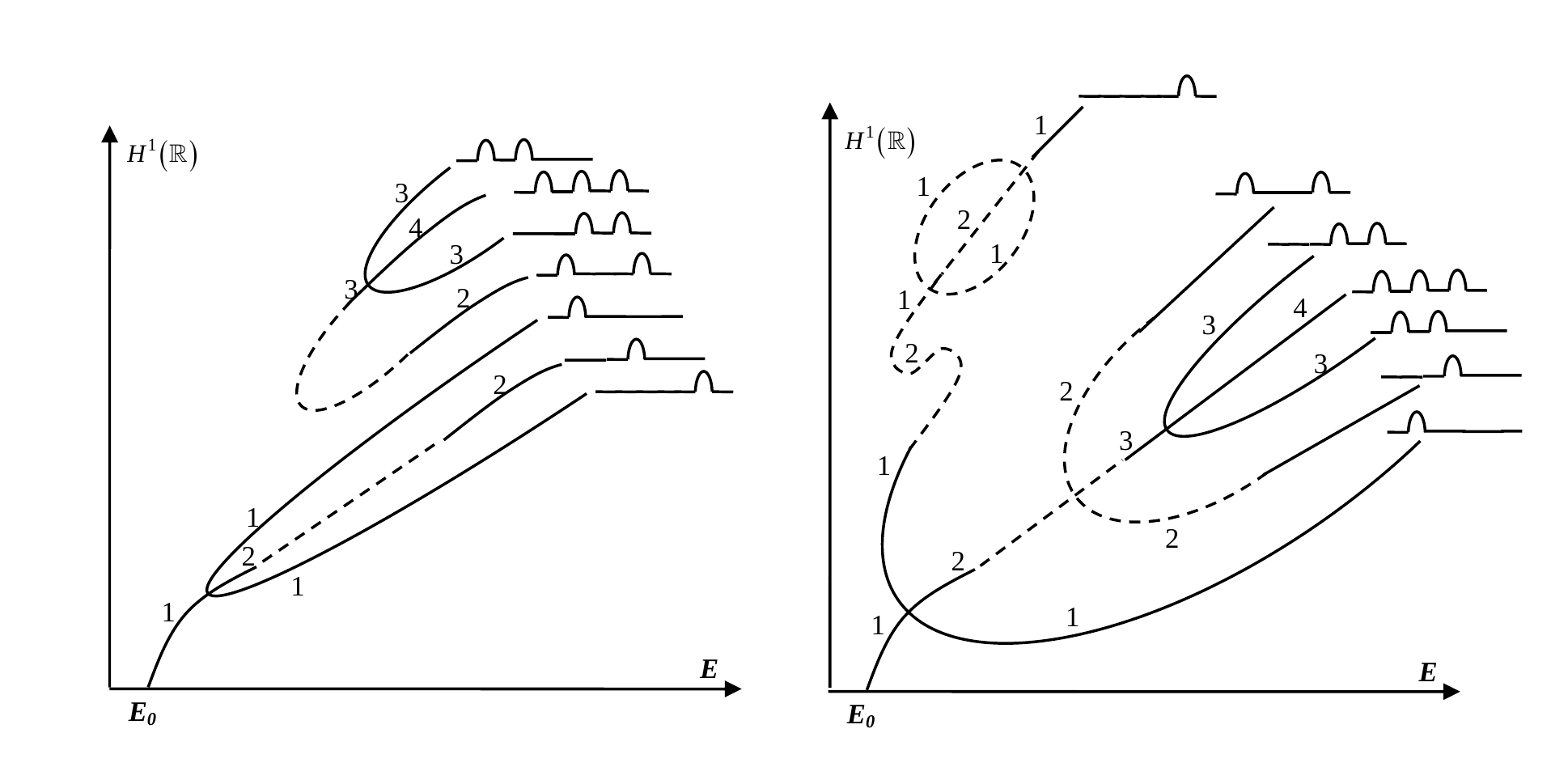}
\end{center}
\caption{Top panel shows two possible ways the even branches
connect. The third is similar to the left panel. Bottom panel shows how the asymmetric branches will
bifurcate from the symmetric ones in the two cases. In all figures
the dotted lines show region where the branches are not completely
understood, i.e. "snaking" or pitchfork like bifurcation may occur
but the latter must lead to loops, see the top branch in the bottom
right panel.} \label{figconnections}
\end{figure}

To avoid the difficult issue of bifurcations from continuous spectrum let us focus first on the attractive case $\gamma <0.$ Suppose we assume that $D_\phi F(\phi_E,E),\ E>E_0,$ is nonsingular along the ground state branch emerging from $(0,E_0).$ If we can now show that the branch can be uniquely continued for all $E>E_0,$ i.e. it approaches the right boundary of the bifurcation diagram, and we can identify the limit point $\lim_{E\rightarrow\infty}\phi_E$, and if the the linearization at the limit point must have two (or more) negative eigenvalues, then we have a contradiction since the linearization had only one negative eigenvalue near $(0,E_0).$ Hence the existence of a singular point along this branch is guaranteed and the resulting bifurcation can give us new branches of ground states. We have actually sketch a result that not only shows there are ground states for all $E>E_0$ but improves on the result of Theorem \ref{th:sb1} by identifying a symmetry breaking bifurcation:

\begin{theorem}\label{th:sb2} Consider an attractive nonlinearity $\gamma <0,$ and a potential $V$ which is invariant under a finite group of Euclidian symmetries on $\R^n.$ Assume that the action of the group is nontrivial on any critical point of $V(x)$ different from $x=0,$ and assume that $x=0$ is a non-degenerate critical point of $V$ different from a minima. Then the branch of ground states bifurcating from $(0,E_0)$ undergoes a second bifurcation past which the symmetric bound states become orbitally unstable. Moreover, one of the new branches emerging from the bifurcation point is made of asymmetric ground states which are generally orbitally stable.
\end{theorem}
In particular, for a double well potential, the bifurcation is of pitchfork type and the emerging branch is made of ground states which localize in one of the two wells, see Figure \ref{figconnections} bottom panels and reference \cite{kkp}.

The theorem is proven using three intermediate results which are important themselves because they determine all ground states provided a few remaining obstacles are surmounted, see Remarks \ref{rmk:r2}, \ref{rmk:exs} and \ref{rmk:r3} below. The first result is:
\begin{theorem}\label{th:R2}
If a $C^1$ branch of coherent states approaches the top or bottom boundary of the bifurcation diagram, i.e. $E\rightarrow E_*,\ 0<E_*<\infty,\ \|\phi_E\|_{H^1}\rightarrow\infty,$ then we have:
$$Q(\phi_E)\rightarrow\infty,\quad\frac{\|\phi_E\|_{L^{p+2}}^{p+2}}{Q(\phi_E)}\rightarrow 0.$$
If a $C^1$ branch of coherent states approaches the right boundary $E\rightarrow\infty,$ then there exists $b>0$ such that $$\frac{\|\phi_E\|_{L^{p+2}}^{p+2}}{E^{2/p+1-n/2}}\rightarrow b,\ \frac{Q(\phi_E)}{E^{2/p-n/2}}\rightarrow -\gamma\frac{(2-n)p+4}{2p+4}b,\ \frac{\|\nabla\phi_E\|_{L^2}^2}{E^{2/p+1-n/2}}\rightarrow -\gamma\frac{npb}{2p+4}.$$
\end{theorem}
These estimates are obtained from the ordinary differential equation valid along these branches:
\begin{equation}\label{eq:denergy}
\frac{d{\cal E}}{dE}(\phi_E)=-E\frac{dQ}{dE}(\phi_E)
\end{equation}
combined with the equation \eqref{eq:ev} and Pohozaev's identity (essentially the $L^2$-scalar product between \eqref{eq:ev} and  $x\cdot\nabla\phi_E$), which leads to closed differential inequalities for $\|\phi_E\|_{L^{p+2}}^{p+2},$ see \cite{kkp} and \cite{kn:bif} for details.

\begin{remark}\label{rmk:r2} The caveat is that the theorem does not yet cover the case of branches undergoing infinitely many bifurcations in all neighborhoods of the boundary (hence they cannot be parametrized by a $C^1$ map in $E$ in any neighborhood of the boundary). However, most of these peculiar situations have been resolved in the sense that they lead to similar estimates which can be used in the next results, see \cite{kn:bif}.
\end{remark}

The theorem is essential in finding the limit points of the bound state branches at the boundary of the bifurcation diagram. For example, at the $\{E\rightarrow E_*,\ 0<E_*<\infty,\ \|\phi_E\|_{H^1}\rightarrow\infty\}$ part of the boundary, the estimate above, combined with the \eqref{eq:tF} form of the equation and the fact that $(-\Delta+1)^{-1}:H^{-1}\mapsto H^1$ is an isomorphism imply that $\phi_E/\sqrt{Q(\phi_E)}$ converges in $H^1$ to a solution of $-\Delta\psi+E_*\psi =0.$ The latter has only the zero solution which contradicts $\|\phi_E/\sqrt{Q(\phi_E)}\|_{L^2}\equiv 1.$

A contradiction is also obtained at the $E=0$ from hypothesis (SA) and comparison principles for the linearized (self-adjoint) operator. The comparison principle relies heavily on the fact that the nonlinearity is always negative, see \cite{kn:bif}.

At the $E\rightarrow\infty$ portion of the boundary the estimates imply that the change of variables:
\begin{equation}\label{chvar}\psi_E(x)=E^{-1/p}\phi_E(E^{-1/2}x+x_0)),\ x_0\in\R^n\end{equation} leads to a uniformly bounded curve $E\mapsto\psi_E$ in $H^1$ and transforms \eqref{eq:ev} into:
 $$-\Delta\psi_E +E^{-1}V(E^{-1/2}x+x_0)\psi_E+\psi_E+\gamma |\psi_E|^{p}\psi_E=0$$
which formally converges to
 \begin{equation}\label{ev0}-\Delta\psi +\psi+\gamma |\psi|^{p}\psi=0.\end{equation}
The rigorous result is:
\begin{theorem}\label{th:R3}
There are no coherent states approaching $\{E=0\}$ and $\{E>0,\ \|\cdot\|_{H^1}\rightarrow\infty\}$ boundary of the bifurcation diagram. Ground states approaching $E\rightarrow\infty$ boundary converge in $H^1$, modulo the re-scaling \eqref{chvar} and rotations in the complex plane, to a superposition of positive solutions of \eqref{ev0} each localized at a critical point of the potential $V.$
\end{theorem}
Note that the result at $E\rightarrow\infty$ has been conjectured in \cite{RW:88}. Our convergence argument uses concentration compactness \cite[Section 1.7]{caz:bk} combined with a rather delicate analysis of bifurcations from infinity, see \cite{kn:bif} for details. For example, if splitting of profiles would occur then at least one of them must move towards infinity, and since $\lim_{|x|\rightarrow\infty}V(x)=0$ we can show that the profile converges to a solution of the equation without potential. But we are dealing with ground states so this solution must be positive modulo rotations. It is known that, modulo translations, there is only one such solution and the properties of the linearized operator at this solution are also known. Using a Lyapunov-Schmidt decomposition based on the linearized operator we show that there are no bifurcations from solutions of the translation invariant problem concentrated at infinity into a solution of our problem with potential under mild hypotheses on the behavior of the potential at infinity.

\begin{remark}\label{rmk:exs} New solutions of translation invariant NLS equation \eqref{ev0} may be discovered based on Theorem \ref{th:R3}. Indeed, the corresponding result for excited states (i.e., solutions of \eqref{eq:ev} which are not ground states) is that as $E\rightarrow\infty$ the re-scaled $\psi_E$ either converges strongly to a superposition of positive solutions of \eqref{ev0}, some of them multiplied by $-1,$ and each localized at a critical point of the potential $V,$ or \eqref{ev0} must have solution that cannot be obtained from the positive one via translations or rotations in the complex plane. There are no such solutions in space dimension $n=1$ (hence the theorem applies to all coherent states in one space dimension) but their existence/non-existence in higher dimensions is an open problem. Note that, in principle, the re-scaled $\psi_E$ can be numerically traced along excited state branches at large $E.$ If profiles that change sign emerge (instead of profiles in which the positive part drifts away from the negative part) then the profile is a new solution of \eqref{ev0}. The algorithm can start from excited states of \eqref{eq:ev} which bifurcate from zero at the second and higher eigenvalues of the linear operator $-\Delta+V.$  The existence of such eigenvalues is guaranteed for sufficiently negative potentials.
\end{remark}

To obtain all limit points of the ground state branches at $E\rightarrow\infty$ we combine Theorem \ref{th:R3} with the local bifurcation result:
\begin{theorem}\label{th:R3'}
At $E=\infty,$ from any superposition of positive or negative solutions of \eqref{ev0} each localized at distinct, non-degenerate, critical points of the potential $V$
bifurcates, modulo the re-scaling \eqref{chvar} and rotations in the complex plane, exactly one curve of coherent states for \eqref{eq:ev}. These coherent states have as many nodal points as the number of sign changes in the superposition. The number of negative eigenvalues of the linearization calculated at these coherent states can be computed with the formula: $k+n_1+n_2+\cdots +n_k$ where $k$ is the number of profiles and $n_j,\ j=1,2,\ldots k$ is the number of negative directions for the Hessian of the potential calculated at the critical point where the $j^{th}$ profile localizes.
\end{theorem}
See Figure \ref{figbranches} for an illustration of this theorem in the case of a double well potential.
The theorem is reminiscent of the result in \cite{oh}, see also \cite{fw,ABC}, for the semiclassical limit. Note that we are not in the semiclassical
limit, our re-scaled equation immediately below \eqref{chvar} differ
by the $E^{-1}$ factor in front of the potential making it an even more degenerate problem. Moreover, our method can be adapted to the semiclassical case and gives a stronger result by not only showing uniqueness of such solutions but also providing their parametrizations and spectral properties of the linearized operator. The extension of Theorem \ref{th:R3'} to degenerate critical points is still open.

\begin{remark}\label{rmk:r3} As of now Theorems \ref{th:R3}-\ref{th:R3'} do not exclude or treat the case of multiple profiles localizing at the {\em same} critical point of $V,$ see \cite{ny} for a related result. For ground states the phenomenon does not occur at local minima, but the other cases are still open.
\end{remark}

The compactness argument at $E\rightarrow\infty$ can be extended in the interior of the bifurcation diagram domain to obtain:
\begin{theorem}\label{th:R4}
Any set of  ground-states $(\phi_E,E)$ which is bounded in  $H^1(\mathbb{R}^n)\times (E_1,\infty),$ where $E_1>0,$ is relatively compact and any limit point is a solution of \eqref{eq:ev}.
 \end{theorem}

Now, Theorem \ref{th:sb2} follows from a contradiction argument. Suppose that along the symmetric branch starting at $(0,E_0)$ no eigenvalues of the linearized operator cross zero. Then, by Lemma \ref{th:R1} and the implicit function theorem the branch can be continued and remains symmetric until it reaches the boundary of the bifurcation diagram. By Theorem \ref{th:R3} it will have $E\rightarrow\infty$ and, in this limit, it will converge, modulo re-scaling \eqref{chvar}, to a superposition of positive solutions of \eqref{ev0} each localized at a critical point of V (some may localize at the same critical point). If the limit is localized at $x=0$ then from Theorem \ref{th:R3'} we deduce that the linearized operator along this branch at large $E$ has at least two negative eigenvalues (one plus the number of negative directions for the Hessian of the potential at $x=0$)in contradiction with the fact that it had only one negative eigenvalue near $E_0.$ If the limit has a profile localizing at a non-zero critical point $x_0,$ by symmetry it must have a profile (positive solution of \eqref{ev0}) at each point in the orbit of $x_0$ under the action of the Euclidian group. In this case the number of negative eigenvalues of the linearized operator at large $E$ is at least the the number of profiles, see Theorem \ref{th:R3'}, which is at least the number of points in the orbit. By hypothesis the latter is at least 2 and gives a contradiction. Consequently, there must be an $E_*,\ E_0<E_*<\infty ,$ such that an eigenvalue of $D_\phi F(\phi_E,E)$ converges to zero as $E\nearrow E_*.$ By Theorem \ref{th:R4} there is a limit point $(\phi_{E_*},E_*)$ and local bifurcation theory can be used to analyze the branches emerging from this point.

More importantly, for analytic nonlinearities ($p$ an even, positive integer), our Theorem \ref{th:R4} combined with global bifurcation theory imply that ground states not only organize themselves
in smooth manifolds but the manifolds can also be smoothly continued past
their singularity points (i.e. bifurcation points) until they either form loops or reach the
boundary of the bifurcation diagram region $H^1(\mathbb{R}^n)\times (0,\infty),$ see \cite{bt,kn:bif}.
Note that if we somehow exclude the cases described in Remarks \ref{rmk:r2} and \ref{rmk:r3} then Theorems \ref{th:R3}-\ref{th:R3'} give us all the limit points at the
boundary and we can now trace back all ground-states.

For example consider the symmetric double well potential in one
space dimension which has three critical points, see top panel in
Figure \ref{figbranches}. We claim that all ground states of this problem are given by the bottom left panel in Figure \ref{figconnections}.
A similar picture can be obtained for excited states with a fixed number of nodal points (zeroes).
Indeed, by first restricting the analysis to the
Banach subspace of even functions in $H^1$ we get via Theorem \ref{th:R3'} (and modulo rotations)
the three curves near $E=\infty$ in addition to the one given by (SA) near $E_0,$
see the left panel of Figure \ref{figbranches}. Global bifurcation theory says
that the latter
connects smoothly with one of the former. Hence, we have three
possibilities, two are presented in the upper panels of Figure
\ref{figconnections}, the third is similar to the left panel. The remaining two curves of symmetric ground states must connect with each other since, again by Theorem \ref{th:R3}, they cannot end up on top or left boundaries of the bifurcation diagram, neither can they end at $(0,E_0)$ due to the uniqueness of the branch emerging from this point. Sturm-Liouville theory allows only a simple eigenvalue to cross zero at each bifurcation point, so, to match the number of negative eigenvalue on the symmetric branches we find that, in the case described in the left panels of Figure \ref{figconnections}, one more bifurcation point is needed on each (the turning point is already a bifurcation on the top curve). Since each of these bifurcations correspond to an antisymmetric eigenvector in the kernel of the linearized operator, one asymmetric branch emerges from each, see bottom left panel in Figure \ref{figconnections}. They already match the curves of asymmetric ground states given by Theorem \ref{th:R3'} at $E=\infty,$ hence the picture is complete. A similar counting argument can be done for the right panels in Figure \ref{figconnections}. Moreover, the same analysis can now be performed for excited states with one nodal point (one change of sign), two nodal points, etc, since Theorem \ref{th:R3} applies to them in one dimension, see Remark \ref{rmk:exs}.

The above analysis did not include the multi-profile ground states which, as $E\rightarrow\infty,$ may have more than one profile localizing at the local maxima of the potential, see Remark \ref{rmk:r3}. Recent numerical investigations in \cite{kkn:bif} show that they are present and the branch starting from one profile at each minima when $E=\infty$ turns back and connects to a branch with two profiles both localizing at the local maxima $x=0,$ while the branch starting from a profile at each critical point connects to the branch with three profiles all of them localizing at $x=0.$  A rigorous understanding of this phenomena is underway.

\begin{remark}\label{rmk:hd} Theorems \ref{th:R2}-\ref{th:R4} are valid in any space dimension, however, to obtain all ground states, the counting argument described above needs to be adapted when non-simple eigenvalues of the linearized operator cross zero. In practical applications the multiplicity of
these eigenvalue is due the Euclidian symmetries of the underlying phenomenon hence its Hamiltonian. The symmetries can be used to simplify the normal form of the local bifurcation, see \cite{goss}. A case by case study is underway, beginning with potentials invariant under finite group of symmetries (such as under reflection w.r.t. hyperplanes, or generated by rotations with a fixed angle) and finishing with potentials invariant under continuous group of symmetries such as spherical ones.
\end{remark}
\begin{remark}\label{rmk:sa} Bifurcations from continuous spectrum may occur in the absence of  the spectral hypothesis {\em (SA)}. Indeed, the sketch of proof for Theorem \ref{th:R3} showed that {\em (SA)} is essential in excluding branches which approach the $\{E=0\}$ boundary. While we can build the picture of all ground states starting now from the branches given by Theorem \ref{th:R3'} at large $E,$ some of these branches will end up at $\{E=0\}.$ To complete the picture we need to find all limit points on this boundary, in particular we need to understand bifurcations from the edge of the continuous spectrum, see Lemma \ref{th:R1}. A summary of recent progress in such problems can be found in \cite{miw:15}. Repelling nonlinearities $\gamma>0$ also fall in this category as preliminary calculations show that all branches of coherent states end up at the left boundary $\{E=0\},$ see \cite{jls:gsd} for a different method applicable to ground-states only. More complicated nonlinearities may also push the coherent states towards this boundary. For example $-|\phi|^2\phi+|\phi|^4\phi$ formally behave like an attractive nonlinearity near $(0,E_0)$ but at large bound-states the repelling part dominates. Hence a turning point is formed on the branch starting at $(0,E_0)$ and the conjecture is that it eventually approaches $\{E=0\},$ see \cite{jac}. The study of such nonlinearities and more general ones is in progress.
\end{remark}
\begin{remark}\label{rmk:pot} Confining potentials $\lim_{|x|\rightarrow\infty}V(x)=\infty$ allow for stronger results compared to Theorems \ref{th:R2}-\ref{th:R4}. Indeed, the bound states now belong to the Banach space $\{\phi\in H^1\ :\ \int_{\R^n}V(x)|\phi(x)|^2dx<\infty\}$ which embeds compactly in $L^2,$ see the previous subsection. This implies that the linearized operator has purely discrete spectrum, that the set of solutions of \eqref{eq:ev} is relatively compact, and that the map $\tilde F,$ see \eqref{eq:tF} is compact. In particular Theorems \ref{th:R2}-\ref{th:R4} are valid for all coherent states. Based on this observation a rigorous study vortices in rotating but confined Bose-Einstein Condensates is underway, see \cite{kgf} for a recent summary of open problems, results and applications.
\end{remark}
\begin{remark}\label{rmk:gbc1} Non-analytical nonlinearities require compactness results stronger than the one given by Theorem \ref{th:R4} i.e., valid also for approximate solutions of \eqref{eq:ev}. Such compactness is needed to construct a degree for the map $\tilde F$ in \eqref{eq:tF} on which the global bifurcation theory relies, see \cite{rab}. Such results hold for confining potentials, see remark above, or repelling nonlinearities $\gamma >0,$ see \cite{jls:bsd,jls:gsd}. The problem for non confining potentials combined with attractive nonlinearities is open.\end{remark}

 \subsection{Orbital Stability}\label{sse:os}

Two of the most cited results in orbital stability of coherent structures are the ones by Grillakis, Shatah and Strauss in \cite{gss,gss1}. One of its refinements \cite{Gr}, which is applicable to the Schr\"{o}dinger and Klein-Gordon equations because of the diagonal structure of their linearization, implies that, in the example presented in the previous subsection, all branches with more than one negative eigenvalue in the spectrum of the linearized operator are unstable while the ones with exactly one negative eigenvalue are stable provided their $L^2$ norm is strictly increasing in $E,$ see figure \ref{figconnections}. However, neither the results in \cite{gss,gss1} nor their numerous refinements cover all possible cases. For example, in the Schr\"{o}dinger case with attractive nonlinearity, the first excited state bifurcating from zero at the second eigenvalue of $-\Delta +V$ is outside the scope of the current orbital stability theory. Thanks to hundreds of pages of proofs based on asymptotic stability techniques, see \cite{ty2} and \cite{swc}, we now know that this branch is unstable in the weakly nonlinear regime provided a resonance condition is satisfied. Is there a simpler way to study the stability of such coherent states, one that will not rely on weak nonlinearities and resonance conditions?

In the general framework presented in Section \ref{se:gnw} the theory uses the Lyapunov functional: $u\mapsto\ {\cal E}(u)-\omega Q(u)$ to study the stability of the coherent states $(\phi_\omega,\omega)$ which are solutions of \eqref{cseq} hence critical points of the Lyapunov functional. The results in \cite{gss,gss1} exploit the fact that $Q$ is invariant under the dynamics and can be summarized as follows: if $\phi_\omega$ is a local minimizer of the Lyapunov functional {\em restricted} to the manifold $Q(u)=Q(\phi_\omega)$ then $\phi_\omega$ is orbitally stable; if $\phi_\omega$ is a saddle point of the Lyapunov functional {\em restricted} to the manifold $Q(u)=Q(\phi_\omega)$ with an odd number of negative directions i.e., odd number of negative eigenvalues of the Hessian {\em restricted} to the tangent space of the manifold $Q(u)=Q(\phi_\omega)$ at $\phi_\omega,$ then $\phi_\omega$ is linearly and orbitally unstable. More precisely, for stability the Hessian can be nonnegative or it can have one negative eigenvalue over the entire space $X$ which disappears when the domain is restricted to the tangent space of the codimension one manifold $Q(u)=Q(\phi_\omega)$ which turns out to be equivalent with the condition $\partial_\omega Q(\phi_\omega)<0$ at the particular $\omega$ under study. However, the theory leaves open the cases when the Hessian has more than one negative eigenvalue over the entire space $X$ and an even number of them remain when restricting to the tangent space of $Q(u)=Q(\phi_\omega).$ The example in the above paragraph is in the case with two negative eigenvalues and no recent refinements of the theory can cover it. Also note that none of the refinements applies to the general framework described in section \ref{se:gnw} but to rather particular cases.

Is it possible to show that if the coherent state $\phi_\omega $ is a saddle point of the Lyapunov functional {\em restricted} to the manifold $Q(u)=Q(\phi_\omega)$ then it is orbitally unstable? Note that the Lyapunov functional is actually invariant under the dynamics, hence initial data on the manifold with energy just below the energy of the coherent state evolve on a level set that takes it far away from the coherent state. More precisely, there is a fixed neighborhood of the orbit of the coherent state which is left by all orbits with initial data approaching the coherent state from a negative direction of the Hessian. This idea has been partially exploited in \cite{gss} but there the negative direction turns out to be an unstable direction of the linearized dynamics $\frac{dv}{dt}=JD^2\left({\cal E}-\omega Q\right)(\phi_\omega)[v]$ i.e., an eigenfunction of a positive eigenvalue of $JD^2\left( {\cal E}-\omega Q\right)(\phi_\omega),$ hence an exponential growth of the distance between orbits leads to instability. This is not the case in the example discussed in the first paragraph of this subsection and in many others. But the point is that even in the absence of unstable directions for the linearized dynamics, the presence of negative direction for the Hessian suffices to prove a (much weaker) linear growth of distance between certain orbits in a small neighborhood of the coherent state which still implies instability. This work is still in progress.

Note that, if the answer to the above question is affirmative then the theory becomes a characterization of orbital stability i.e., in the general framework of \eqref{hpde} that is only invariant under the action $T(s)$ of a {\em one dimensional Lie group}, $s\in\mathbb{R}$, the coherent state $\phi_\omega$ given by \eqref{cseq} is orbitally stable if and only if it is a local minimizer of the Lyapunov functional: $u\mapsto\ {\cal E}(u)-\omega Q(u)$ {\em restricted} to the manifold $Q(u)=Q(\phi_\omega)$ i.e., the Hessian over the whole space $X$ of the Lyapunov functional can have at most one negative eigenvalue (counting multiplicity) and if it has one then $\partial_\omega Q(\phi_\omega)$ must be negative at the particular $\omega$ under study. Since the Hessian is basically the linearization of the equation for coherent structures \eqref{cseq}, its eigenvalues change continuously along manifolds of solutions and cross zero only at bifurcation points. Therefore, the stability properties can now be deduced directly form the bifurcation diagrams, see for example figure \ref{figconnections}.

\section{Asymptotic Stability of Coherent States}\label{se:as}

When the techniques proposed in the previous Section lead to a new branch of orbitally stable coherent states for \eqref{hpde} or a bifurcation point involving both stable and unstable branches, the question is whether the dynamics of solutions
starting near the branch can be described in detail. In
particular, asymptotic stability would mean that the solutions
converge to certain coherent structures on the branch but in a weaker norm corresponding to a space $Z,\ X\hookrightarrow Z,$ see the discussion at the end of Section \ref{se:gnw}. The methods
to uncover the convergence are dynamical in nature and, near a branch of coherent structures and {\em away} from bifurcation points, can be
summarized as follows: one decomposes the solution into a finite
dimensional evolution on the manifold of coherent structures and a
correction
  \begin{equation}\label{nlspdec}
u(t)=T(\omega(t)t+\theta(t))[\phi_{\omega(t)}+u_d(t)]
\end{equation}
where the parameters $\omega(t),\ \theta(t)\in\R$ are to be chosen later.
Then the equation for the correction becomes:
\begin{equation}\label{eq:ud1}
\frac{du_d}{dt}=JL_{\omega(t)}u_d+G(\omega(t),u_d(t))
\end{equation}
where $L_{\omega}$ is the linearization (with respect to $u$) of $D_u {\cal E}(u)-\omega D_u Q(u)$ at
$\phi_\omega $ and $G$ contains only quadratic and higher order
terms in $u_d.$ Most of the times equation
\eqref{eq:ud1} is be analyzed via a Duhamel formula using the
propagator of a {\em fixed } linearization i.e., $W(t)u_0$ solves:
$$\frac{du}{dt}=JL_{\omega_0}u,\qquad u(0)=u_0$$
and
\begin{equation}\label{duh}u_d(t)=W(t)u_d(0)+\int_0^tW(t-s)[JL_{\omega(s)}-JL_{\omega_0}]u_d(s)ds+\int_0^tW(t-s)G(\omega(s),u_d(s))ds.\end{equation}
In this case $\omega(t),\ \theta(t)$ are chosen such that $u_d(t)$ is always in the invariant subspace of $JL_{\omega_0}$ that complements its (two dimensional) null space. In the absence of other eigenvalues the invariant space corresponds to the continuous spectrum, and the advantage is that, on this subspace,
estimates of type:
\begin{equation}\label{dest}
\|W(t)\|_{Z^*\mapsto Z}\sim |t|^{-\alpha},\qquad \alpha >0
\end{equation}
where $Z$ is a Banach space with $X\hookrightarrow Z$ densely, were
already available. The disadvantage is the presence of the linear
term in the second integral. In fact, for the particular case of
NLS, the linear term lead to restrictions on the nonlinearity to
supercritical regimes $p>4/n$ in \eqref{eq:nse}, see \cite{sw:mc1,sw:mc2,PW,bp1,bp2,bp3,cuc,sz},
or, when Stricharz type estimates were used, to
critical and supercritical regimes $p\ge 4/n$, see \cite{gnt,miz1,miz2}.
The results in \cite{kz1,kz2,km1,km2} show that, for small solitary
waves in NLS, from estimates \eqref{dest} one can obtain estimates
of the same type for the propagator of the {\em time dependent}
linear operator in \eqref{eq:ud1}. Hence one can use:
\begin{equation}\label{duh1}u_d(t)=\tilde W(t,0)u_d(0)+\int_0^t\tilde W(t,s)G(\omega(s),u_d(s))ds.\end{equation}
where $\tilde W(t,s)u_0$ solves the non-autonomous equation
\begin{equation}\label{lpr}\frac{du}{dt}=JL_{\omega(t)}u,\qquad u(s)=u_0\end{equation}
and now the parameters $\omega(t),\ theta(t)$ are chosen such that $u_d(t)$ is always in the invariant subspace of $JL_{\omega(t)}$ that complements its null space, hence, in the absence of other eigenvalues, on this subspace we have:
\begin{equation}\label{dest1}
\|\tilde W(t,s)\|_{Z^*\mapsto Z}\sim |t-s|^{-\alpha},\qquad \alpha >0
\end{equation}
As a result the restriction to critical and supercritical nonlinearities has been lifted.

Essential in this approach is to obtain estimates \eqref{dest1} from
\eqref{dest}. While technical in nature this step can be generalized
to other equations, see \cite{bk}, because it relies only on
$V(t)=JL_{\omega(t)}-JL_{\omega_0}$ being a small, localized in
space (but time dependent and maybe complex valued) scatterer (potential) and on
strong dispersion of the linearized equation, i.e.
\begin{equation}\label{alpha}\alpha\ge 1.\end{equation}
Smallness is not necessary since orbital stability implies that
large deviation in $V(t)$ are only along the orbits of the coherent
structures which can be mod out, see \cite{kz3}. Space localization will
always be present when dealing with solitary waves i.e., the scatterer is a power of the solitary wave. The only
hypothesis that cannot yet be relaxed is \eqref{alpha}, in
particular the method is inapplicable in one dimensional NLS.



A much more delicate dynamics occurs near an intersection of stable and unstable manifolds of coherent structures. Subsection \ref{sse:bm} shows that existence of such bifurcation
points is generally the rule rather than the exception, hence
understanding the dynamics around bifurcations is a necessary step
in studying asymptotic completeness. Note that finite time behavior of small
solitary waves near the bifurcation point discovered first in \cite{kksw} has been studied for example in
\cite{MW10} via an approximation with a finite dimensional dynamical system.

Unfortunately, the recent progress in asymptotic stability near an orbitally stable coherent state cannot yet determine the full dynamical picture near a bifurcation point. Current results, see \cite{bc,cum,zg}, rely on a spectral assumption for the linearization $JL_\omega$ which fails at the bifurcation point. More precisely, as one approaches the bifurcation point $(\omega\rightarrow\omega_*)$ from a stable branch, two, purely imaginary and complex conjugate eigenvalues (which can be non-simple) of $JL_\omega$ approach zero (this corresponds to one eigenvalue, maybe non-simple, of $L_\omega=D^2{\cal E}(\phi_\omega)-\omega D^2Q(\phi_\omega)$ approaching zero, see the discussion in section \ref{sse:bm}). Past the bifurcation point they move from zero back up on the imaginary axis (along on the stable branch) or into a positive and negative eigenvalue (along the unstable branch). The above cited results may be applicable for some $\omega\not=\omega_*$ along the stable branch but not to all. This is because the two eigenvalue which approach zero as $\omega\rightarrow\omega_*$ have multiples close to any other eigenvalue, therefore violating the discrete spectrum non-resonance condition required by current results. G. Zhou has done yet unpublished work that may remove the non-resonance condition. Even if this work is vetted the results only say that for each $\omega\not=\omega_*$ along the {\em stable branch} there is a small ball in the Hilbert space $X$ centered at $\phi_\omega$ such that initial data from the ball asymptotically converge (in the weaker norm) to a coherent state (not necessarily $\phi_\omega$). However, the radius of this ball is related to the distance $d$ between the smallest eigenvalue (in absolute value) and zero, and goes to zero as $\omega$ approaches the bifurcation point because of the approx. $1/d$ changes of variable necessary to bring the system to a normal form that uncovers the {\em radiation damping mechanism} which leads to the decay of the projection of the solution onto the invariant subspace of this smallest eigenvalue via a resonant interaction  with the radiative part corresponding to the projection onto the continuous spectrum. In conclusion, asymptotic stability can be shown only in a conical neighborhood of the stable branch with vertex at the bifurcation point (and the vertex excluded).

What happens with initial data in a ball centered at the bifurcation point which, of course, has relatively large regions not contained in the conical neighborhoods described above? A first step would be to determine the stable invariant manifolds along the unstable branch. This is obtained via implicit function theorem type results from the invariant subspace of the linearization that complements the one corresponding to the positive eigenvalue. Initial data on this manifold will converge to an unstable coherent state, see \cite{ty2} for a related result. Outside this manifold one can use a spectral decomposition of the dynamics with respect to the linearization at the bifurcation point $JL_{\omega_*}.$ Note that the invariant subspace corresponding to the zero eigenvalue contains the kernel of $D^2{\cal E}(\phi_\omega)-\omega D^2Q(\phi_\omega)$ which caused the bifurcation to occur in the first place, moreover the projection on this kernel parameterizes the branches near the bifurcation point via the standard Lyapunov-Schmidt decomposition. Since the initial data is away from the stable invariant manifold corresponding to the unstable branch one expects a short time exponential growth of the projection onto the direction of the stable branch. Once this projection becomes dominant a change of variables can be employed in order to use the linearization and associated spectral decomposition at the (time dependent) stable coherent state given (parameterized) by the values of this dominant projection. In this coordinates the techniques discussed in the above paragraph  are expected to lead to asymptotic convergence towards a stable coherent state. This is work in progress.

\section{Conclusions}

While variational methods are inappropriate to study {\em all coherent states} of a given nonlinear wave equation \eqref{hpde} recent progress in NLS equation shows promise for bifurcation methods. We learned from the example described in subsection \ref{sse:bm} that one can start from any known solution of \eqref{cseq} and expect that it is part of a smooth manifold of solutions that can be extended to the boundary of the domain (subset of $X\times\mathbb{R})$) inside which  the linearization of \eqref{cseq}, $D^2{\cal E}(u)-\omega D^2Q(u),$ is Fredholm (a finite dimensional kernel suffices as this is a self-adjoint operator). Such results are the main theorems of the current global bifurcation theory as developed by Rabinowitz, Dancer, Toland, Buffoni and others. However, this is not enough to discover bifurcation points along the manifold or new branches of solutions. But if one has spectral information about the linearization at the starting point and at the end point of the manifold (both on the boundary) then existence of bifurcations can be proven and the branches emerging from them can be studied. Ideally one would want to find all limit points of manifolds of coherent states at the boundary of the Fredholm domain. Now, one can start from any such limit point and use the fact that the manifold approaching it can be continued (via global bifurcation theory) until it reaches another limit point on the same boundary (or the same point in case a loop forms). Loops can sometimes be ruled out via bifurcation in cones type arguments, see \cite{bt}, and the symmetry or spectral properties of the initial limit point can severely reduce the choices of the end point such that a numerical investigation or a rigorous theorem can determine it as in the example. Once the two end points of the manifold are determined the change in number of negative eigenvalues of the linearized operator at the two ends can tell us the number of eigenvalues crossing zero hence the number and type of bifurcations along the manifold. These bifurcations can now be analyzed via Lyapunov-Schmidt decompositions and normal forms, see for example the singularity theory in \cite{goss}. The emerging branches of solutions have also limits at the boundary of the Fredholm domain. The process repeats until all branches of solutions are found.

There are four essential steps in the method summarized above:
  \begin{itemize}
  \item[(R1)] identify the domain in the $X\times\mathbb{R}$ space where the
  linearized operator of \eqref{cseq} is Fredholm;
  \item[(R2)] inside this domain show relative compactness of either the set of solutions of \eqref{cseq} or of a map for which the set of fixed points coincides with the set of solutions of \eqref{cseq};
  \item[(R3)] identify all limit points of solution branches on the boundary of the Fredholm domain;
  \item[(R4)] find the rules by which the limit points connect
  via manifolds of solutions inside the domain and characterize the bifurcation points along these manifolds.
  \end{itemize}

(R1) is already done for most wave equations as linearization
is used in any analysis or numerical simulation for nonlinear
problems. While this method focuses on finding coherent states inside the domain where the linearization is Fredholm, note that outside it one can sometimes show non-existence of coherent states with certain localization properties based on smoothness of the spectral measure of the linearized operator, see \cite{bl}.

(R2) is a technical step but an essential assumption in global bifurcation theory. It also helps with the (R3) step. For problems with real analytic energy, relative compactness of the solution set of \eqref{cseq} suffices and implies that the only obstacle for unique continuation of any manifold of coherent states (via implicit function theorem) is that it either reached the boundary of the Fredholm domain or a bifurcation (singularity) point. In the latter case, the manifold reemerges on the other side of the bifurcation point because of the structure of zeroes of analytical maps. Compactness and structure of zeroes for analytical maps combine again to prevent the existence of infinitely many singularities in bounded domains. Hence the manifold either reaches the boundary of the Fredholm domain or forms a loop, see \cite{bt}. If the energy is non-analytical a stronger form of compactness is required. Equation \eqref{cseq} is transformed into a fixed point problem, for example $\phi=(-\Delta + 1)^{-1}\psi,\ \psi=(1-E)\phi-V(x)\phi-\gamma
|\phi|^{p}\phi\stackrel{def}{=}\tilde F(\psi, E)$ in the NLS case \eqref{eq:ev}, and the map $\tilde F:X^*\mapsto X^*$ (or from $Y$ to $Y$ where $X\hookrightarrow Y\hookrightarrow X^*$) is required to be relatively compact, see \cite{rab}.

Note that the stronger type of compactness comes for free when $X\hookrightarrow Y$ is compact, which is the case for waves on bounded domains, or when the range of $\tilde F$ is made of functions with a prescribed decay at infinity, which is the case for problems with confining potentials, $\lim_{|x|\rightarrow\infty}V(x)=\infty,$ or with repelling nonlinearities see \cite{jls:bsd,jls:gsd}. For NLS with attractive nonlinearities, a delicate argument involving concentration compactness and non-existence of bifurcations from profiles concentrated at infinity is used in \cite{kn:bif} to obtain relative compactness of the set of solutions. It may be adapted for general wave problems since $X, Y$ are usually Sobolev spaces.

(R3) can be split into two parts: first obtain rigorous estimates for coherent structures approaching
the boundary of the Fredholm  domain and use compactness arguments to identify possible limit points, then use local
bifurcation theory from these limit points to identify all nearby branches of solutions. For the first part one can use the energy and charge $Q$ together with the identity $\frac{d{\cal E}}{d\omega}(\phi_\omega)=\omega\frac{dQ}{d\omega}(\phi_\omega)$ and the equation \eqref{cseq} both valid for coherent states. The goal is to obtain closed differential inequalities for the terms in the energy which can lead estimates in certain limits i.e., as the branch approaches different parts of the boundary of the Fredholm domain. In the NLS example with attractive nonlinearity, the quadratic terms in the energy \eqref{Enls} grow much faster than the superquadratic term corresponding to the nonlinearity provided the $H^1$ norm blows up and  $\omega$ remains finite. Equation \eqref{eq:ev} becomes linear in this limit and it turns out that the limiting equation has no solution, hence the absence of coherent states near this boundary. At the $\omega\rightarrow -\infty$ boundary the kinetic and nonlinear terms grow faster than the potential term which leads to the limiting equation \eqref{ev0} via the re-scaling \eqref{chvar}. The limiting points are then among the solutions of the limiting equation. In the example we were looking at positive coherent structures (ground states) and the limiting equation has exactly one such solution modulo translations, hence the limiting points described in figure \ref{figbranches}. In general, one can focus on examples for which the limiting equations are well studied. However if the procedure leads to poorly understood equations this theory is a motivating factor in studying them.

The second part i.e., identify the nearby branches from the limit points, may be solved via standard local bifurcation theory when the linearization at the limiting points (which are solutions of the limiting equation) is Fredholm. However, this is not the case when coherent states approach the part of the boundary where the linearization has zero at an edge of the essential spectrum. This happens for repelling nonlinearities, see \cite{jls:bsd,jls:gsd}, for attractive nonlinearities when the linearization of \eqref{cseq} at $u=0$ has no discrete spectrum and especially for Dirac equation, regardless of nonlinearity, since the linear Dirac operator has essential spectrum everywhere except a bounded interval. In all these cases local
bifurcations from the edge of the essential spectrum must be understood
in order to find all limit points on this boundary and complete the bifurcation diagram. These are
notoriously difficult problems but promising results in this direction are described in \cite{miw:15}.

(R4) amounts to grouping the limiting points based on which closed subspace of $X$ they belong to, usually based on symmetry properties such as the subspace of even functions in our example. Since global bifurcation theory applies in any Banach space each limit points must connect to another one in the same group. If there are more than two in a subgroup then not only numerical simulations can help but also rigorous arguments arguments combining the mismatch in the negative eigenvalues of the linearization at the two endpoints with the type of bifurcations supported by the eigenspaces corresponding to the eigenvalues that cross zero as one moves from one limiting point to the other. Note that problems invariant under finite and continuous groups of Euclidian symmetries have non-simple eigenvalues in the spectrum of the linearized operator. If they cross zero a classification of bifurcations induced by them is necessary before we can proceed to identify all coherent states.

The bifurcation method relies and provides information on the spectrum of the linearized operator along the manifolds of coherent states. There is already a rich theory that uses the spectral information to determine the stability of coherent states and the long time dynamics of nearby solutions. While a resolution of the {\em Asymptotic Completeness Conjecture}, see Section \ref{se:intro}, still seems far away, it appears that a systematic study of all coherent states supported by nonlinear wave equations, their bifurcation points, their stability and the nearby dynamics is within reach.

\end{document}